\documentclass[11pt ]{cocv}
\usepackage{latexsym,amssymb,amsmath}
   \usepackage{times}
   \usepackage{float}
\newif\ifpdf
\ifx\pdfoutput\undefined
\pdffalse 
\else
\pdfoutput=1 
\pdftrue
\fi
\ifpdf
\usepackage[pdftex]{graphicx}
\else
\usepackage[dvips]{graphicx}
\fi
\newcommand{\N} {\mathbb{N}}             %
\newcommand{\R} {\mathbb{R}}             %
\newcommand{\iR} {\int_{\R} }      
  \newcommand{\preuve}{\noindent\textit{ Proof -~}}
 \newcommand{\findemo}{\hfill $\Box$}
  
 \newcommand{\bv}{ BV(\Omega)}
 \newcommand{\W}{ \Omega}
 \newcommand{\ldeux}{L^2(\Omega)}

\newtheorem{theo}{Theorem}[section]
\newtheorem{prop}{Proposition}[section]

\newtheorem{rem}{Remark}[section]

\begin{document}
\ifpdf
\DeclareGraphicsExtensions{.pdf,.jpg}
\else
\DeclareGraphicsExtensions{.eps,.jpg}
\fi

\title{An active curve approach for tomographic reconstruction of
  binary radially symmetric objects}

\author{I. Abraham}\address{ 
CEA Ile de France-
BP 12\\
91680 Bruy\`eres le Ch\^atel\\
FRANCE\\
isabelle.abraham@cea.fr}

\author{R. Abraham} \author{M. Bergounioux} \address{Laboratoire MAPMO-
F\'ed\'eration Denis Poisson\\
Universit\'e d'Orl\'eans\\
BP 6759\\
ORLEANS cedex 02\\
FRANCE\\
romain.abraham@univ-orleans.fr, ~maitine.bergounioux@univ-orleans.fr}
 \today

\begin{abstract}
 This paper deals with a method of tomographic reconstruction of
 radially symmetric objects from a
 single radiograph, in order to study the behavior of shocked
 material. The usual tomographic reconstruction algorithms such as
 generalized inverse or filtered back-projection cannot be applied
 here because data are very noisy and the inverse problem
 associated to single view tomographic reconstruction is highly
 unstable. In order to improve the reconstruction, we propose here to
 add some a priori assumptions  on the looked after object. One of
 these assumptions is that the object is binary and consequently, the
 object may be described by the curves that separate the two
 materials. We present a model that lives in BV space and leads to a non local Hamilton-Jacobi equation, via a level set strategy.  Numerical experiments are performed (using level
 sets methods) on synthetic objects.
 \end{abstract}
 
         
   \maketitle

\section{Introduction}   

Medical scanner is the most used application of tomographic
reconstruction. It allows to explore the interior of a human body. In
the same way, industrial tomography explores the interior of an object
and is often used for non-destructive testing. 

We are interested here in a very specific application of tomographic reconstruction
for a physical experiment described later. The goal of this experiment
is to study the behavior of a material under a shock. We obtain
during the deformation of the object an X-ray radiography by high
speed  image capture. We suppose
  this object is radially symmetric, so that one
radiograph  is enough to reconstruct the 3D object.

\begin{figure}[H]
\begin{center}
   \includegraphics[width=12cm]{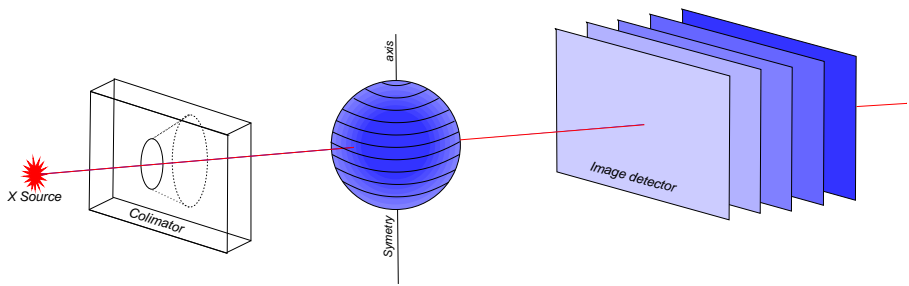} 
   \caption{Experimental setup}
   \end{center}
   \end{figure}

Several authors have proposed techniques (standard  in medical tomography) for
tomographic reconstruction when enough projections (from different
points of view) are available: this allows to get an analytic formula for the solution
(see for instance \cite{Her} or \cite{Dus}). These methods cannot be used
directly when only a few number of projections is known. Some
alternative methods have been proposed in order to partially
reconstruct the densities (see for instance \cite{Din}).
We are interesting here in single view tomographic
reconstruction for radially symmetric object (see for instance
\cite{Han} for a more complete presentation of the subject). As any tomographic
reconstruction, this problem leads to an ill-posed inverse problem. As
we only have one radiograph, data are not very redundant and
the ill-posed character is even more accurate.

We present here a tomographic method adapted to this specific
problem, originally developed in
\cite{AA}, and based on a curve evolution approach. The main idea is
to add some a priori knowledge on the object we are studying in order
to improve the reconstruction. The object may then be described by a
small set of characters (in this case, they will be curves) which are
estimated by the minimization of an energy functional. This work is
very close to another work by Feng and al \cite{FKC}. The main
difference is the purpose of the work: whereas they are seeking
recovering textures, we are looking for accurate edges.
It is also close to the results of Bruandet and al \cite{BPDB}.
However, the present work handles very noisy
data and highly unstable inverse problems, and shows how this method is
powerful despite these perturbations. Further, we take here into
account the effects of blur (which may be non-linear) and try
to deconvolve the image during the reconstruction.

Let us mention at this point that our framework if completely different
from the usual tomographic point of view, and usual techniques (such
as filtered back-projection) are not adapted to our case. Indeed,
usually, as the X-rays are supposed to be parallel (this is also the
case here), the ``horizontal'' slices of the object are supposed to be
independent and studied separately. Usual regularization techniques
deal with one slice and regularize this particular slice. Here,
because of the radial symmetry, the slices are composed of concentric
annulus and do not need any regularization. The goal of this work is
to add some consistency between the slices in order to improve the 
reconstruction.

The paper is organized as follows. First we present the
physical experiment whose data are extracted and explain what are the
motivations of the work. Next, we introduce the projection operator. In Section 4, we present a continuous model with the suitable functional framework and prove existence result.  Section 5 is devoted to   formal computation  of the energy derivative in order to state some optimality conditions. In Section 6,  a front propagation point of view is adopted  and   the level set method leads to a non local  Hamilton-Jacobi equation.
In the last section, we present some numerical results and give hints for numerical schemes improvement. 

\section{Experiment}

This work is part of some physical experiments whose goal is the
study of the behavior of shocked material. The present experiment consists in making a
hull of well known material implode using surrounding explosives. The
whole initial physical setup (the hull, the explosives ...) are radially
symmetric. A reasonable assumption is to suppose that during the
implosion, everything remains radially symmetric.

Physicists are looking for the shape of the interior at some fixed
time of interest. At that time, the interior may be composed of
several holes which also may be very irregular. Figure
\ref{fig:object} is a synthetic object that contains all the
standard difficulties that may appear. These difficulties are
characterized by:
\begin{itemize}
\item Several disconnected holes.
\item A small hole located on the symmetry axis (which is the area
  where the details are difficult to recover).
\item Smaller and smaller details on the boundary of the top hole in
  order to determine a lower bound detection.
\end{itemize}

To achieve this goal, a X-rays radiograph is obtained. In order to
extract the desired informations, a tomographic reconstruction must be
performed. Let us note here that, as the object is radially symmetric,
a single radiography is enough to compute the reconstruction.

A radiography measures the attenuation of X-rays through the object. A
point on the radiography will be determined by its coordinates $(u,v)$
in a Cartesian coordinates system where the $v$-axis will be the
projection of the symmetry axis. If
$I_0$ is the intensity of the incident X-rays flux, the measured flux
$I$ at a point $(u,v)$ is
given by
$$I=I_0 e^{-\int\mu(r,\theta,z)d\ell}$$
where the integral operates along the ray that reaches the point
$(u,v)$ of the detector, $d\ell$ is the infinitesimal element of
length along the ray and $\mu$ is the local attenuation coefficient.
For simplicity, we will consider that
this coefficient is proportional to the material density. To deal with
linear operators, we take the Neperian logarithm of this attenuation
and will call the transformation
$$\rho\longmapsto \int\rho d\ell$$
the projection operator.

Through the rest of the paper, in order to simplify
the expression of the projection operator, we will suppose that the
X-ray source is far enough away from the object so that we may
consider that the rays are parallel, and orthogonal to the symmetry axis. As a consequence, the horizontal
slices of the object may be considered separately to perform the projection.
 
\begin{figure}[H]
\centerline{\includegraphics[width=5cm,angle=270]{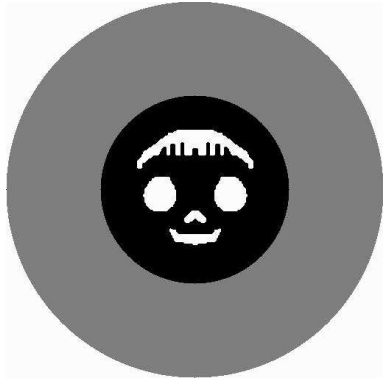}}

\vspace{0.5cm}

\caption{Slice of a typical binary radially symmetric object by a
  plane that contains the symmetry axis (the $z$-axis).
}\label{fig:whole-object}
\end{figure}

\begin{figure}[H]
\centerline{\includegraphics[width=5cm,angle=270]{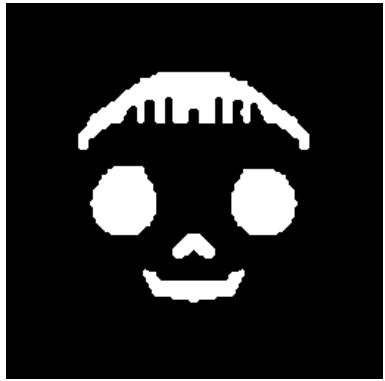}}

\vspace{0.5cm}

\caption{Zoom on the interior of the object of Figure~\ref{fig:whole-object}.
The homogeneous material is in black whereas the holes are in white.}\label{fig:object}
\end{figure}

As the studied object is radially symmetric, we will work in a system
of cylinder coordinates $(r,\theta,z)$ where the $z$-axis is the
symmetry axis. The object is then described by the density at the
point $(r,\theta,z)$, which is given by a function $f$ which depends
only on $(r,z)$ by symmetry. In the text, the notation $f$ will always
refer to the density of the object. A typical function $f$ is given in
Figure \ref{fig:whole-object}. It represents an object composed of concentric
shells of homogeneous materials (called the ``exterior'' in what
follows) surrounding a ball (called the ``interior'') of another homogeneous material that
contains some empty holes. This figure may be viewed as a slice of the
object by a plane that contains the symmetry axis. To recover the
3D-object, it suffices to perform a rotation of this image around the
$z$ axis. For instance, the two round white holes in the center are in
fact the slice of a torus. As the looked-after characteristic of the object is the
shape of the holes, we will focus only on the interior of the object
(see Figure \ref{fig:object}). We here handle only binary objects
composed of one homogeneous material (in black) and some holes (in white).

\section{A variational approach}
\subsection{The projection operator}
We first explicit the projection operator and its adjoint. 
\begin{prop}
In the case of a radially symmetric object, the projection operator,
denoted by $H$, is given, for  every function $f \in L^\infty(\R_+\times\R)$ with compact support, by  
\begin{equation}\label{H}
\forall (u,v)\in \R\times \R \qquad Hf(u,v)=2\int_{|u|}^{+\infty}f(r,v)\frac{r}{\sqrt{r^2-u^2}}dr.
\end{equation}
\end{prop}
\preuve Consider a 3D-object which is described by a function $\tilde f(x,y,z)$ (Cartesian coordinates system). The projection operator $H$ is 
$$ H\tilde f(u,v) =\iR \tilde f(x, u,v) \, dx~.$$
In the case of radially symmetric object, we parametrize the object by a function $f(r,z)$ with cylinder coordinates. Therefore
$$ \tilde f(x,y,z) = f(\sqrt{x^2+y^2},z) ~.$$
Then, we have, for $u\ge 0$
$$ Hf(u,v) = \iR \tilde f(x, u,v) \, dx =\iR f(\sqrt{x^2+u^2},v) \, dx = 2 \int_0^{+\infty} 
f(\sqrt{x^2+u^2},v) \, dx ~.$$
We perform  the following change of variable 
$$ r= \sqrt{x^2+u^2},~x \ge 0 \Longleftrightarrow x= \sqrt{r^2 - u^2}, ~r \ge u~,$$
to get
$$  Hf(u,v) = 2 \int_u^{+\infty}  f(r,v) \,\frac{r}{ \sqrt{r^2 - u^2}} \,dr ~.$$
For $u < 0$, we have 
$$  Hf(u,v) = 2 \int^u_{-\infty}  f(r,v) \,\frac{|r|}{ \sqrt{r^2 - u^2}} \,dr ~.$$
Using the change of variable $u :to -u$ and the fact that $r \mapsto f(r,v)$ is even,  by symmetry, we get
$$Hf(u,v)=2\int_{|u|}^{+\infty}f(r,v)\frac{r}{\sqrt{r^2-u^2}}dr.$$
\findemo
\begin{rem} Operator $H$ may be defined by density on measurable functions $f$
  such that all the partial applications $f(\cdot,z)$ belong to
  $L^2(\R_+)$. Then, all the functions $Hf(\cdot,v)$ belong to the space $BMO(\R)$ of bounded mean oscillations functions :
$$ \mbox{BMO}(\R) = \{ f :\R \to \R~|~Ê\sup_{R>0} \left( \frac{1}{R} \int_{|x-y|< R} \left | f(x) - f_R (y)\right | \, dx \right ) < +\infty~\}~,$$
where $\displaystyle{f_R(y) = \frac{1}{2R} \int_{y-R}^{y+R} f(x) \, dx }$. For more details, one can refer to \cite{Stein}.
\end{rem}

In the sequel, we will need to handle functions $f$ that are defined
on $\R^2$ (instead of on $\R_+\times\R$). We thus define the operator
$H$ for function $f\in L^\infty(\R^2)$ with compact support by
$$Hf(u,v)=2\int_{|u|}^{+\infty}f(sgn(u)r,v)\frac{r}{\sqrt{r^2-u^2}}dr$$
although this has no more physical meaning. Here, the function $sgn$
is defined by
$$sgn(x)=\begin{cases}
1 & \mbox{if }x\ge 0,\\
-1 & \mbox{if } x<0.
\end{cases}$$

We shall also need the back-projection that is the adjoint operator $H^*$ of $H$; it can be computed in a similar way. 
\begin{prop} \label{propadj}The  adjoint operator (in $L^2$) $H^*$ of the  projection operator is given,   for every function $g \in L^\infty(\R^2)$ with compact support by :  
\begin{equation}\label{Hadj}
\forall  r \in ”\R~,~\forall z \in \R, \qquad H^*g(r,z)=2\int_{0}^{|r|}g(sgn(r)u,z)\frac{|r|}{\sqrt{r^2-u^2}}du.
\end{equation}
\end{prop}
\preuve The adjoint operator $H^*$ of $H$ is the unique
operator such that, for every $f$ and $g$ in $L^\infty(\R^2)$ with compact support,
$$\int_{-\infty}^{+\infty} \int_{-\infty}^{+\infty} Hf(u,v) \,g(u,v)  \, dv \, du =  \int_{-\infty}^{+\infty} \int_{-\infty}^{+\infty}  f(r,z) \, H^*g(r,z)  \, dr \, dz ~.$$
Using (\ref{H}) and Fubini's theorem, we get 
\begin{align*}
\int_{-\infty}^{+\infty} \int_{-\infty}^{+\infty} & Hf(u,v) \,g(u,v)  \, dv \, du\\
& = 2 \int_{u=-\infty}^{u=0} \int_{v=-\infty}^{v=+\infty} \int_{r=|u|}^{+\infty}
f(-r,v)\frac{r}{\sqrt{r^2-u^2}} \,g(u,v)  \, dv \, du \, dr \\
& \qquad\qquad +2\int_{u=0}^{u=+\infty}\int_{v=-\infty}^{v=+\infty} \int_{r=u}^{+\infty}
f(r,v)\frac{r}{\sqrt{r^2-u^2}} \,g(u,v)  \, dv \, du \, dr \\
&
=2\int_{v=-\infty}^{v=+\infty}\int_{r=0}^{r=+\infty}\int_{u=-r}^{u=0}\frac{r}{\sqrt{r^2-u^2}}f(-r,v)g(u,v)\,du\,dr\,dv\\
& \qquad\qquad
+2\int_{v=-\infty}^{v=+\infty}\int_{r=0}^{r=+\infty}\int_{u=0}^{u=r}\frac{r}{\sqrt{r^2-u^2}}f(r,v)g(u,v)\,du\,dr\,dv\\
&
=2\int_{v=-\infty}^{v=+\infty}\int_{r=0}^{r=-\infty}\int_{u=r}^{u=0}\frac{-r}{\sqrt{r^2-u^2}}f(r,v)g(u,v)\,du\,(-dr)\,dv\\
& \qquad\qquad
+2\int_{v=-\infty}^{v=+\infty}\int_{r=0}^{r=+\infty}\int_{u=0}^{u=r}\frac{r}{\sqrt{r^2-u^2}}f(r,v)g(u,v)\,du\,dr\,dv\\
&
=2\int_{v=-\infty}^{v=+\infty}\int_{r=-\infty}^{r=+\infty}\int_{u=0}^{u=r}\frac{r}{\sqrt{r^2-u^2}}f(r,v)g(u,v)\,
du\,dr\,dv\\
&=2\int_{v=-\infty}^{v=+\infty}\int_{r=-\infty}^{r=+\infty}\int_{u=0}^{u=|r|}\frac{|r|}{\sqrt{r^2-u^2}}f(r,v)g(sgn(r)u,v)\,
du\,dr\,dv
\end{align*}
So we obtain the following expression for the back projection : 
$$ H^*g(r,z)  =2\int_{0}^{|r|}\frac{|r|}{\sqrt{r^2-u^2}} \,g(sgn(r)u,z)  \, du   ~.$$
\findemo 

\subsection{Toward a  continuous model}

Thanks to the symmetry, this operator characterizes the Radon
transform of the object and so is invertible; one radiograph is enough
to reconstruct the object. The inverse operator is given, for an
almost everywhere differentiable function $g$ with compact support, and for every $r>0$, by
$$H^{-1}g(r,z)=-\frac{1}{\pi}\int_r^{+\infty}\frac{\frac{\partial
  }{\partial x}g(x,z)}{\sqrt{x^2-r^2}}\,dx.$$

Because of the derivative term, the operator $H^{-1}$ is not
continuous. Consequently, a small variation on the measure $g$ leads
to significant errors on the reconstruction. As our radiographs are
strongly perturbed, applying $H^{-1}$ to our data leads to a poor
reconstruction. Due to the experimental setup they are also two main perturbations:
\begin{itemize}
\item A blur, due to the detector response and the X-ray source spot
  size.
\item A noise.
\end{itemize}
Others perturbations such as scattered field, motion blur... also exist
but are neglected in this study.
\\
We denote by $F$ the effect of blurs. We will consider the
following   simplified case where $F$ is supposed to be
linear
\begin{equation}\label{blur}
F(k)=N*k\end{equation}
where $*$ is the usual convolution operation, $k$ is the projected image  and $N$ is a positive
symmetric kernel. 

\begin{rem}
A more realistic case stands when the convolution
operates on the intensity :  then $F$ is of the form
$$F(k)=\frac{-1}{\nu}\ln \left(e^{-\nu k}*N\right)$$
where $\nu$ is the multiplicative coefficient between the density and
the attenuation coefficient. Some specific
experiments have been carried out to measure the blur
effect. Consequently, we will suppose
 that, in both cases, the kernel $N$ is known. The linear blur is not
 realistic but is treated here to make the computations simpler for
 the presentation. 
 \end{rem}

\begin{figure}[ht]
\begin{center}
\includegraphics[width=5cm,angle=270]{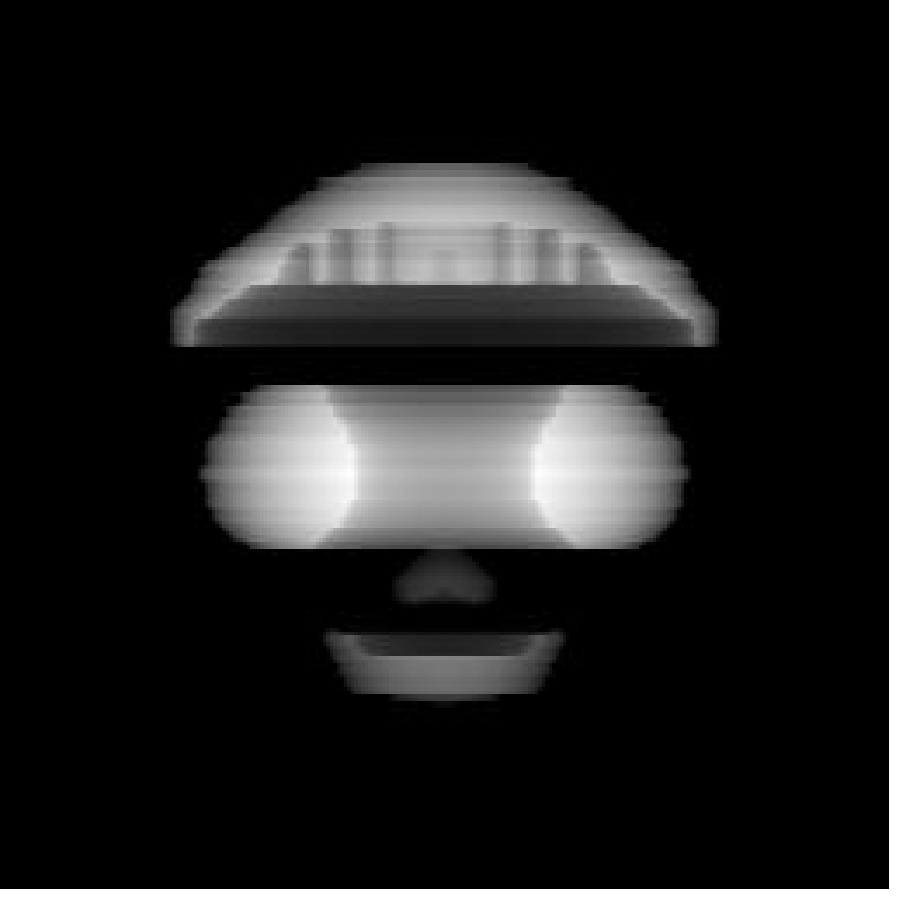}\qquad
\includegraphics[width=5cm,angle=270]{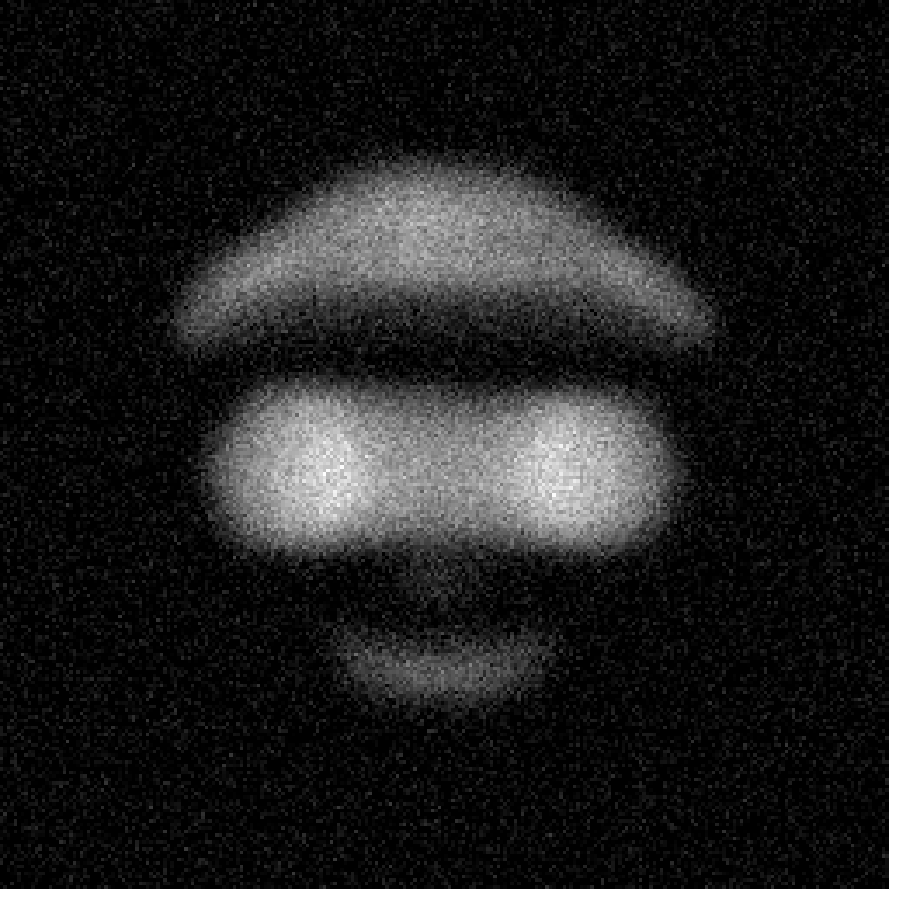}
\caption{Left-hand side: theoretical projection $F(Hf)$ of the object of
  Figure~\ref{fig:object}. Right-hand side: real projection of the
  same object with
  realistic noise and blur.}\label{fig:real_proj}
\end{center}
\end{figure}

\begin{figure}[ht]
\begin{center}
\includegraphics[width=5cm,angle=270]{tete2.eps}\qquad
\includegraphics[width=5cm,angle=270]{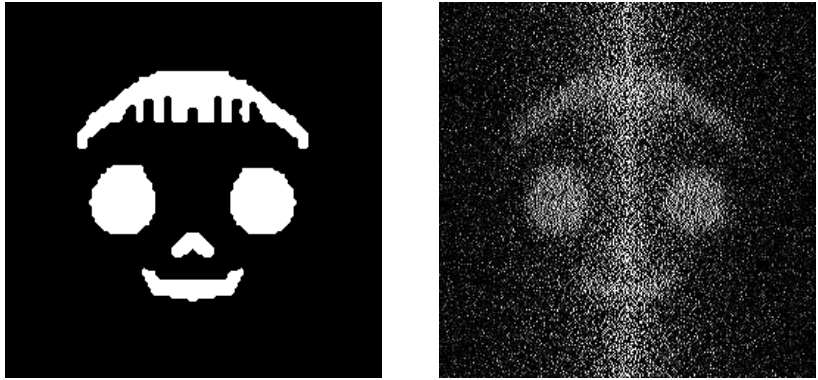}
\caption{Comparison between the real object on the left-hand side and
  the reconstruction computed with $H^{-1}$ applied to the real
  projection on the right-hand side.}\label{fig:inverse}
\end{center}
\end{figure}

The noise is supposed for simplicity to be an additive Gaussian white
noise of mean 0, denoted by $\varepsilon$.
Consequently, the projection of the object $f$ will be
$$g=F(Hf)+\varepsilon.$$
The comparison between the theoretical projection $F(Hf)$ and the
perturbed one is shown on Figure \ref{fig:real_proj}. The
reconstruction using the inverse operator $H^{-1}$ applied to $g$ is
given by Figure \ref{fig:inverse}. The purpose of the experiment is to
separate the material from the empty holes and consequently to
precisely determine the frontier between the two areas, which is
difficult to perform on the reconstruction of Figure \ref{fig:inverse}.

It is clear from Figure \ref{fig:inverse} that the use of the inverse
operator is not suitable. In order to improve the reconstruction, we
must add some a priori knowledge on the object to be
reconstructed. Indeed the object that we reconstruct must satisfy some
physical property.

We chose to stress on two points:
\begin{itemize}
\item The center of the object is composed of one homogeneous known material's density
  with some holes inside.
\item There cannot be any material inside a hole.
\end{itemize}
In a previous work \cite{Lag}, J.M. Lagrange reconstructed the
exterior of the object. In this reconstruction, the density of the
material at the center of the object is known, only the holes are not
reconstructed. In other words, we can reconstruct an object without
holes and we can compute (as $H$ and $F$ are known) the theoretical
projection of this reconstruction. We then act as if the blurred
projection was linear and subtract the projection of the non-holed
object to the data.
In what follows, we will call
experimental data this subtracted image which corresponds to the
blurred projection of a ``fictive'' object of density 0 with some holes
of known ``density'' $\lambda > 0$. Consequently, the space of admitted objects will
be the set of functions $f$ that take values in $\{0,\lambda\}$. This
space of functions will be denoted $\mathcal{F}$ in the sequel.

The second hypothesis is more difficult to take into account. We chose
in this work to tackle the problem via an energy minimization method
where the energy functional is composed of two terms: the first one
is a matching term , the second one is a penalization
term which tries to handle the second assumption. The
matching term will be a $L^2$-norm which is justified by the
Gaussian white noise.

\begin{rem} In the case where $F$ is not linear, the exact method to
remove the exterior is to operate the blur function on the addition of
the known exterior and the center. For the sake of simplicity, we will not
use this method and will consider that the errors are negligible when
subtracting the projections.
\end{rem}
 
  Let us first describe more
precisely the set $\mathcal{F}$. This time, the functions
$f\in\mathcal{F}$ will be defined on $\R^2$, with values still in
$\{0, \lambda\}$, with compact support. Therefore, such a function $f$
  will be characterized by the knowledge of the curves that limit the
  two areas where $f$ is equal to $\lambda$ and to $0$. Indeed, as the
  support of the function $f$ is bounded, these curves are disjoint
  Jordan curves and the density of the inside is $\lambda$ whereas
  the density of the outside is 0. Consequently, the energy that we
  will consider will be a function of $\gamma_f$ where $\gamma_f$ is a set
  of disjoint Jordan curves. For mathematical reasons, we must add an
  extra-assumption : the curves $\gamma_f$ are $\mathcal{C}^1$ so that
  the normal vector of the curves is well-defined (as an orthogonal vector to the tangent one).

In this continuous framework, the matching term is just
the usual $L^2$-norm between $Hf$ and the data $g$ (where $H$ is given by (\ref{H})). So, the first term is
$$E_1(\gamma_f)=\|{F(Hf)-g}\|_2^2.$$
For the penalization term,  we choose
$$E_2(\gamma_f)=\ell(\gamma_f)$$
where $\ell(\gamma_f)$ denotes the length of the curves $\gamma_f$. Let us
remark that this penalization term may be also viewed as the total
variation (up to a multiplicative constant) of the function $f$ because of the binarity.
Eventually, the total energy functional is
\begin{equation}\label{energie}
E(\gamma_f)= \|{F(Hf)-g}\|_2^2+\alpha\ell(\gamma_f)
\end{equation}
which is an adaptation of the well-known Mumford-Shah energy
functional introduced in \cite{MS}. The ``optimal'' value of $\alpha$
may depend on the data.

\subsection{A continuous model in BV space}

The previous analysis gives the mains ideas for the modelization. Now, we make it precise using an appropriate functional framework. 
 Let  $\Omega$ be a bounded  open subset $\R^2$ with Lipschitz  boundary.  We  shall consider bounded variation functions.  Recall that the space of such functions is 
$$\bv = \{ u \in L^1(\Omega) ~|~J(u) < + \infty\,\} $$
where
\begin{equation}\label{ju}
J(u) = \sup \left \{ \int_\Omega u(x) \, \mbox{div }\xi (x) \, dx  ~|~\xi \in \mathcal{C}^1_c(\Omega)~,~\|\xi\|_\infty\le 1~ \right\},
\end{equation}
where $\mathcal{C}^1_c(\Omega)$ denotes the space of $\mathcal{C}^1$ functions with compact support in $\Omega$. The space $\bv$ endowed with the norm 
$$ \|u\|_{\bv} = \|u\|_{L^1} + J(u)~,$$
is a Banach space.
\\
If $u\in \bv$ its derivative in $\mathcal{D}'(\W) $ (distributions) is a  bounded Radon measure denoted $Du$  and  $J(u)$ is the total variation of  $|Du|$ on $\Omega$. Let us recall useful properties of BV-functions 
(\cite{ambrosio}):

\begin{prop}  Let  $\Omega$ be an open subset $\R^2$ with Lipschitz  boundary. \\
1.  If $u\in  \bv$, we get the following decomposition for $Du$ :
$$ Du = \nabla u \,dx + D^su~,$$
where $\nabla u dx $ is the absolutely continuous part of $Du$ with respect of the  Lebesgue measure and  $D^s u$ is the singular part.\\
2. The map $u \mapsto J(u)$ from  $\bv$ to $\R^+$ is lower semi-continuous (lsc) for the $L^1(\Omega)$ topology.\\
3. $\bv \subset L^2(\Omega)$ with compact embedding.\\
4. $\bv \subset L^1(\Omega)$  with compact embedding.
\end{prop}
  
We precise hereafter an important continuity property of the projection  operator $H$. 
\begin{prop}\label{H1} The projection  operator $H$ is continuous from  $L^{2+s}(\Omega)$  to  $L^p(\Omega)$ for every $ p \in [1,+\infty]$ and  $s >0$.
\end{prop}
\preuve It is a direct consequence of H\"older inequality. Let be
$s >0$ and $ f \in L^{2+s}(\W)$. Its  support  is included in  $\bar
\W$ which is included in some $[-M,+M] \times [-M,+M] $ where $M >0$
and only depends on $\W$. It is clear that  $Hf$ is defined everywhere
on $\W$  and, for every $u\ge 0$
$$ 
\begin{array}{rl}
 |Hf(u,v)| &=\displaystyle{ 2\, \left |\int_{u}^{+\infty}f(r,v)\frac{r}{\sqrt{r^2-u^2}}\,dr \right|}
  =\displaystyle{ 2\left |\int_{u}^{M}f(r,v)\frac{r}{\sqrt{r^2-u^2}}\,  dr  \right|}
  \\[0.5cm]
 &\le 2\displaystyle{ \left[ \int_{u}^{M}|f(r,v)|^{2+s} \, dr \right]^{\frac{1}{2+s}} \, \left[\int_{u}^{M}
 \frac{r^q}{(r +u)^{\frac{q}{2} }}\, \frac{1}{ (r -u)^{\frac{q}{2} }}\, dr\right] ^{\frac{1}{q}} }
\end{array} $$
where $q= 1 + \displaystyle{\frac{1}{1+s}}$.  
Therefore
$$ 
  |Hf(u,v)|  
 \le 2\|f\|_{L^{2+s}}  \,  \left[\int_{u}^{M}\frac{r^q}{ r^{\frac{q}{2} }   (r -u)^{\frac{q}{2} }} dr\right ]  
 \le 2 \,   M^{\frac{q}{2} }\|f\|_{L^{2+s}} \,\left[\int_{u}^{M}   (r -u)^{-\frac{q}{2}} dr \right] 
 $$
The computations in the case $u<0$ are similare and lead to the same
inequality (with some additional absolute values).
As 
$$1- \frac{q}{2} = \frac{s}{2(1+s)} >0~,$$
we get 
\begin{equation}
 |Hf(u,v)|  \le  2\, M^{\frac{q}{2} } \|f\|_{L^{2+s}}  \, \frac{2(1+s)}{s} \, [M-u]^\frac{s}{2(1+s)}\le C(\W,s) \|f\|_{L^{2+s}}~;\end{equation}
here and in the sequel  $C(\W,s)$ denotes a generic constant depending on  $s$ and $\W$.
So
 $$ \|Hf \|_\infty \le C(\W,s) \|f\|_{L^{2+s}}~.$$
 As  $\W$ is bounded, this yields 
 \begin{equation}\label{borneH}
 \forall f  \in L^{2+s}(\W) ,~\forall p \in [1,+ \infty] \quad \|Hf \|_{L^p(\W)}   \le C(\W,s) \|f\|_{L^{2+s}}~.\end{equation}
\findemo
 
 As we consider the length of curves,  the most suitable functional space  to set a variational formulation of the reconstruction problem is   $\bv$.
Therefore, we consider the following minimization problem
$$(\mathcal{P}) \qquad \left \{
\begin{array}{l}
\min \|FH f - g \|_2^2 + \alpha J(f) \\
f \in \bv \\
|f(x) | = 1 \mbox{ a.e. on  } \W \end{array}\right . $$
Here
\begin{itemize}
\item $\|\cdot\|_2$  stands for the $\ldeux$- norm, $g \in \ldeux$ and $\alpha >0.$
\item The operator $F$ is  given by (\ref{blur}).  Without loss of generality,  we may assume (for simplicity) that $F = I.$
\item At last,  ``$|f(x) | = 1 \mbox{ a.e.  sur } \W $'',  is the binarity constraint.  We have mentionned that the image 
  $f_o$ takes its values in  $\{0,  \lambda\}$ where  $\lambda>0$. With the change of  variable  $\displaystyle{f =- \frac{2}{\lambda} f_o  +1 }$, we may assume that the image values belong to  $\{ -1, 1\} $. 
\end{itemize}
\begin{rem}A similar problem has been studied in \cite{CKP}  with smoother projection operator and \textbf{convex} constraints. This is not our case. The pointwise constraint 
  `` $|f(x) | = 1 \mbox{ a.e. on } \W $'' is a very hard constraint. The constraint set is not convex and its interior is empty for most usual topologies.  \end{rem}
  Now we may give the main result of this section : 
\begin{theo} Problem  $(\mathcal{P})$ admits at least a  solution. \end{theo}
\preuve  Let $\varphi_n \in \bv$ be a minimizing sequence. It satisfies $\|\varphi_n\|_\infty =1$ ; so
\begin{equation} \label{estim1}
\forall p \in [1, +\infty[,~\forall n \in \N, \qquad \|\varphi_n \|_{L^p} \le |\W|^{\frac{1}{p}}~.\end{equation}
Therefore the sequence  $(\varphi_n)$ is $L^1(\W)$- bounded. As   $J(\varphi_n)$ is bounded as well, the sequence is bounded in  $\bv$.  Thus it converges (extracting a subsequence) to some $\varphi \in \bv$ for the weak-star topology. 
\\
Estimate  (\ref{estim1}) implies the weak    convergence  of $(\varphi_n)$ to $\varphi$ in  $L^{2+s}(\W)$ for every $s >0$. Thanks to the  $H$  continuity property of   proposition \ref{H1},  we assert that  $H\varphi_n$ weakly converges to  $H\varphi$ in $\ldeux$. We get \begin{equation} \label{estim2}
 \|H\varphi   - g \|_{L^2}^2  \le \liminf_{n \to \infty} \|H\varphi_n  - g \|_{L^2}^2  ~,\end{equation}
with the lower semi-continuity of the norm.
\\
Moreover  $\bv$ is compactly embedded in $L^1(\W)$. This yields that $(\varphi_n)$ strongly converges to $\varphi$ in $L^1(\W)$.  As  $J$ is lsc with respect to  $L^1(\W)$- topology, we get 
\begin{equation} \label{estim3}
J(\varphi )\le \liminf_{n \to \infty}  J(\varphi_n)  ~,\end{equation}
Finally
\begin{multline*}
\inf (\mathcal{P}) = \lim_{n\to +\infty} \|H\varphi_n  - g
\|_{L^2}^2+\alpha J(\varphi_n)   \\
\ge \liminf_{n \to \infty} \|H\varphi_n  - g \|_{L^2}^2+\alpha J(\varphi_n)  \ge \|H\varphi  - g \|_{L^2}^2 + \alpha J(\varphi)~.\end{multline*}
As the pointwise constraint  is obviously satisfied, $\varphi$ is a solution to  $(\mathcal{P})$. 
 \findemo

\section{Computation of the energy derivative }

Now we look for optimality conditions. Unfortunately we cannot compute easily the derivative of the energy in the $\bv$ framework. 
Indeed we need  regular curves  and  we do not know if the $\bv$ minimizer  provides a curve with  the required regularity. Moreover, the set of constraints is not convex and it is not easy to compute the G\^ateaux- derivative (no admissible test functions).

So we have few hope to get classical optimality conditions and we rather compute minimizing sequences. We focus on particular ones that are given via the gradient descent method inspired by \cite{MS}.  Formally, we look for a family of curves $(\gamma_t)_t\ge 0 $ such that 
$$ \frac{\partial E}{\partial \gamma}  (\gamma_t) \leq 0$$
 so that   $ E(\gamma_t)$ decreases as $t \to +\infty$. 
Let us compute the energy variation when
we operate a small deformation on the curves $\gamma$. In other word,
we will compute the G\^ateau derivative of the energy for a small
deformation $\delta \gamma$:
$$\frac{\partial E}{\partial \gamma}  (\gamma) \cdot \delta \gamma  =\lim_{t\to 0}\frac{E(\gamma
  +t\delta\gamma)-E(\gamma)}{t}\cdot$$
 
\begin{figure}[ht]
\begin{center}
 \includegraphics[scale=0.5]{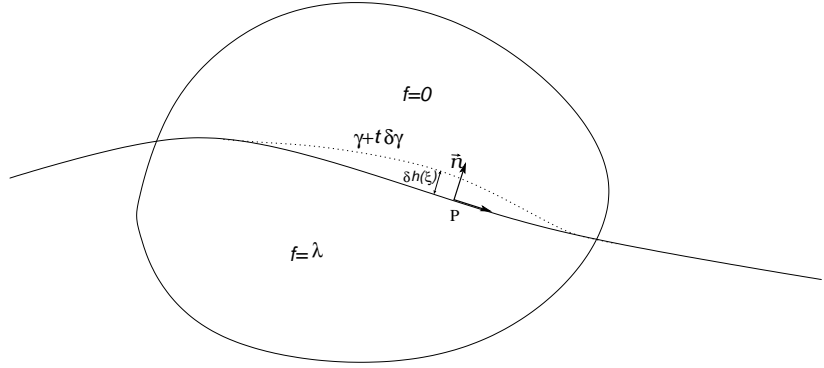}
\caption{Description of a local deformation of the initial curve
  $\gamma$. $P$ is the current point, $U$ is the neighborhood of $P$
  in which the deformation is restricted to and $\gamma+t\delta\gamma$ is
  the new curve after deformation. The interior of the curve is the set where $f=\lambda$}\label{fig:deformation}
\end{center}
\end{figure}

We will first focus on local deformations $\delta\gamma$. Let 
$(r_0,z_0)$ be a point $P$ of $\gamma$. We consider a local  reference
system  which center  is  $P$ and axis are given by the tangent and
normal vectors at $P$ and we denote by $(\xi,\eta)$ the new generic
coordinates in this reference system. With an abuse of notation, we
still denote $f(\xi,\eta)=f(r,z)$. We apply the implicit functions
theorem to parametrize our curve: there exist a neighborhood $U$ of
$P$ and a $\mathcal{C}^1$ function $h$ such that, for every $(\xi,\eta)\in
U$, 
$$(\xi,\eta)\in \gamma\iff \eta=h(\xi) .$$
  Eventually, we get a
neighborhood $U$ of $P$, a neighborhood $I$ of $\xi_0$ and a
$\mathcal{C}^1$ function $h$ such that
$$\gamma\cap U=\Bigl\{(\xi,\eta)\in \R^2\bigm| \eta=h(\xi),\ \xi\in I\Bigr\}.$$
The local parametrization is oriented along  the outward normal $\vec n$ to the curve $\gamma$ at point $P$ (see figure  
\ref{fig:deformation}).  More precisely, we define the local coordinate system $(\vec \tau, \vec n)$ where $ \vec \tau$ is the usual tangent vector, $\vec n$ is the direct orthonormal vector;  we set the curve orientation so that $\vec n$  is the outward normal. The function $f$ if then defined on $U$ by
$$ f(\xi,\eta)=\left \{ \begin{array}{ll}
 \lambda & \mbox{if }\eta <  h(\xi)\\
0 & \mbox{if }\eta\ge h(\xi)
\end{array}\right. $$
This parametrization is described on figure \ref{fig:deformation}. We
then consider a local (limited to $U$)  deformation $\delta\gamma$ along the normal vector. This is equivalent to handling a $\mathcal{C}^1$ function $\delta
h$ whose support is included in $I$. The new curve $\gamma_t$ obtained
after the deformation $t\delta \gamma$ is then parametrized by
$$ 
\eta=\left \{ \begin{array}{ll}
h(\xi)+t\delta h(\xi) & \mbox{for }(\xi,\eta)\in U\\
\gamma & \mbox{otherwise}
\end{array}\right. $$
This defines a new function $f_t$:
\begin{equation}\label{eq:deltaf}
f_t(\xi,\eta)=\left \{ \begin{array}{ll}
 f(\xi,\eta) & \mbox{if }(\xi,\eta)\not \in U\\
 \lambda & \mbox{if }(\xi,\eta)\in U\cap\{\eta< h(\xi)+t\delta h(\xi)\}\\
0 & \mbox{if }(\xi,\eta)\in U\cap\{\eta\ge  h(\xi)+t\delta h(\xi)\}\\
\end{array}\right.
\end{equation}
We will also set $\delta f_t=f_t-f$. This deformation is described on figure \ref{fig:deformation}.
 
The G\^ateau derivative for the energy $E_2$ has already been computed
in \cite{MS} and is
$$\frac{\partial E_2}{\partial \gamma}(\gamma) \delta \gamma=-\int
_\gamma\mbox{curv}(\gamma)\bigl(\xi,h(\xi)\bigr)\delta h(\xi)d\xi$$
where $\mbox{curv}$ denotes the curvature of the curve and $\delta h$ is the parametrization of $\delta \gamma$. 
\\
It remains to compute
the derivative for the matching term. First we estimate  $\delta f_t$: a simple computation shows that
$$ \delta f_t (\xi,\eta) = \left \{ \begin{array}{rl}
0&\mbox{if } \eta \ge h(\xi) + t\delta h (\xi) \mbox{ or } \eta(\xi)\le h(\xi)  \\
\lambda&\mbox{if }h(\xi) \le  \eta \le h(\xi) + t\delta h (\xi) \end{array}\right. \mbox{ in  case } \delta h  \ge 0~ $$
 and
 $$ \delta f_t (\xi,\eta) = \left \{ \begin{array}{rl}
0&\mbox{if } \eta \le h(\xi) + t\delta h (\xi) \mbox{ or } \eta(\xi)\ge h(\xi)  \\
- \lambda&\mbox{if }h(\xi) \ge  \eta \ge h(\xi) + t\delta h (\xi) \end{array}\right. \mbox{ in  case } \delta h  \le 0~ $$
Now we compute 
 $ E_1(\gamma_t) - E_1(\gamma)$ where $\gamma$ (resp. $\gamma_t$) is the curve associated to the function $f$ (resp. $f_t$): 
 
\begin{align*}
E_1(\gamma_t) & - E_1(\gamma)\\
 & = \iR \iR \left ( (g-FHf_t)^2 -(g-FH f)^2 \right) (u,v)\, du \,dv\\
 &=  \iR \iR \left ( (g-FHf-FH\delta f_t)^2 -(g-FH f)^2 \right)
 (u,v)\, du \,dv\\
 & = -2 \iR \iR   (g-FHf)(u,v) \, FH\delta f_t  (u,v)\, du \,dv
 +\underbrace{ \iR \iR (FH \delta f_t)^2(u,v) \, du\, dv}_{= o(t)}~ \\
&= -2\left < (g-FHf),FH\delta f_t\right >_{L^2}+o(t) \\
&=-2\left < H^*F^*(g-FHf),\delta f_t\right >_{L^2}+o(t).
\end{align*}

To simplify the notations, we denote by $\mathcal{A}f:=(H^*F^*g - H^*F^* F
H f)$ so that we need to compute
$$\lim_{t\to 0}\frac{1}{t}\langle \mathcal{A}f,\delta f_t\rangle_{L^2}.$$

As $\delta f_t$ is zero out of the neighbourhood $U$, we have
$$\langle \mathcal{A}f,\delta f_t\rangle_{L^2}=\int
_U\mathcal{A}f(\xi,\eta)\delta f_t(\xi,\eta)d\xi\,d\eta.$$

In the case $\delta h\ge 0$, we have,
$$\langle \mathcal{A}f,\delta f_t\rangle_{L^2}=\lambda\int_{\xi\in
  I}\int_{\eta=h(\xi)}^{\eta=h(\xi)+t\delta h(\xi)}\mathcal{A}f(\xi,\eta)d\xi\,d\eta.$$

As the function $\mathcal{A}f$ is continuous (and thus bounded on
$U$), we may pass to the limit by dominated convergence and get
$$\lim_{t\to 0}\frac{1}{t}\langle \mathcal{A}f,\delta
f_t\rangle_{L^2}=\lambda\int_I\mathcal{A}f\bigl(\xi,h(\xi)\bigr)\delta h(\xi)d\xi.$$

In the case $\delta h<0$, we have
$$\langle \mathcal{A}f,\delta
f_t\rangle_{L^2}=\int(-\lambda)\int_{\xi\in I}\int_{\eta=h(\xi)+t\delta
  h(\xi)}^{\eta=h(\xi)}\mathcal{A}f(\xi,\eta)d\xi\,d\eta$$
and we obtain the same limit as in the nonnegative case.

Finally, the energy derivative is
$$
\frac{\partial E}{\partial \gamma} (\gamma_f)\cdot \delta \gamma_f =-2\lambda \int_I\mathcal{A}f\bigl(\xi,h(\xi)\bigr) \delta
h(\xi)\, d\xi -  \alpha \int_I
\mbox{curv}(\gamma_f) (\xi,h(\xi) ) \delta
h(\xi) \, d\xi
$$
  If we set $\beta = \displaystyle{\frac{\alpha}{2}} $, we get
  \begin{equation}\label{derivEbis}
\frac{\partial E}{\partial \gamma} (\gamma_f) \cdot \delta \gamma_f=-  2\int_I \left ( \lambda \mathcal{A}f + \beta  
\mbox{curv}(\gamma_f) (\xi,h(\xi) )\right )  \delta h(\xi) \, d\xi~. 
\end{equation} 
As  $ \delta h = < \delta \gamma_f , \vec n >$  formula  (\ref{derivEbis}) may be written 
\begin{equation}\label{eq:derivee}
\frac{\partial E}{\partial \gamma} (\gamma_f) \cdot \delta \gamma_f =- 2 \int _\gamma  \left (  \lambda \mathcal{A}f + \beta  
\mbox{curv}(\gamma_f) (s)\right )  < \delta \gamma_f , \vec n > \, ds
\end{equation}
where $\vec n$ denotes the outward pointing normal unit vector of the curve
$\gamma$, $<\cdot,\cdot>$ denotes the usual scalar product in
$\R^2$ and $c(s)$ is a positive coefficient that depends on the
curvilinear abscissa $s$.

The latter expression is linear and
continuous in $\delta \gamma$, this
formula is also true for a non-local deformation (which can be
achieved by summing local deformations).

\section{Front propagation and  level set method}

The goal of the present section is to consider a family of curves
$(\gamma_t)_{t\ge 0}$ that will converge toward a local minimum of
the functional energy. From equation (\ref{eq:derivee}), it is
clear that if the curves $(\gamma_t)$ evolve according to the
differential equation
\begin{equation}\label{eq:edp}
\frac{\partial \gamma}{\partial t}=  (\lambda \mathcal{A}f + \beta  
\mbox{curv}(\gamma_f) )\vec n,
\end{equation}
the total energy will decrease.

To implement a numerical scheme that discretizes equation
(\ref{eq:edp}), it is easier to use a level set method (see \cite{Set}
for a complete exposition of the level set method). Indeed, equation
(\ref{eq:edp}) may present some instabilities, in particular when two
curves collide during the evolution or when a curve must
disappear. All these evolutions are handled easily via the level set method.

The level set method consists in viewing the curves $\gamma$ as the
0-level set of a smooth real function $\phi$ defined on $\R^2$. The
function $f$ that we are seeking is then just given by the formula
$$f(x)= \lambda 1_{\phi(x)>0}.$$
We must then write an evolution PDE for the functions $\phi_t = \phi(t,\cdot)$ that
corresponds to the curves $\gamma_t$. Let $x(t)$ be a point of the curve
$\gamma_t$ and let us follow that point during the evolution. We know
that this point evolves according to equation \ref{eq:edp}
$$x'(t)= \left( \lambda \mathcal{A}f + \beta  
\mbox{ curv}(\gamma_f) \right) (x(t))\vec n.$$
We can re-write this equation in terms of the function $\phi$ recognizing
that
$$\vec n=\frac{\nabla \phi }{|\nabla \phi|}\quad\mbox{and}\quad  \mbox{curv}(\gamma)=\mathrm{div}\left(\frac{\nabla \phi }{|\nabla \phi |}\right)$$
where $\nabla $ stands for the gradient of $\phi$ with respect to $x$, $|\cdot|$ denotes the euclidean norm. The evolution equation
becomes
$$x'(t)= \left( \lambda  \mathcal{A}\left(\lambda 1_{\phi(t, \cdot)> 0}\right)+\beta\mbox{ div}\left(\frac{\nabla 
  \phi }{|\nabla \phi |}\right)\right)\frac{\nabla \phi }{|\nabla 
  \phi |}(t,x(t)).$$
Then, as the point $x(t)$ remains on the curve $\gamma_t$, it satisfies
  $\phi_t(x(t))=\phi(t,x(t))= 0$. By differentiating this expression, we obtain
$$\frac{\partial\phi}{\partial t}+\bigl<\nabla\phi, x'(t)\bigr>=0$$
which leads to the following evolution equation for $\phi$:
$$
\frac{\partial\phi}{\partial t}+ |\nabla_x\phi  |\left(\lambda  \mathcal{A}\left(\lambda 1_{\phi(t,\cdot)>0}\right)+\beta\mathrm{ div}\left(\frac{\nabla 
  \phi}{|\nabla \phi |}\right)\right)= 0~,
$$
that is 
\begin{equation}\label{eq:levelset}
\frac{\partial\phi}{\partial t}=|\nabla \phi |\left( \lambda^2 H^*F^*FH\left(1_{\phi(t,\cdot)>0}\right)-\beta\mathrm{ div}\left(\frac{\nabla 
  \phi}{|\nabla \phi |} \right)- \lambda H^*F^*g\right).
\end{equation}

The above equation is an Hamilton-Jacobi equation which involves a non local term (through $H$ and $F$). Such equations are difficult to handle especially when it is not monotone (which is the case here).  In particular, existence and/or uniqueness of solutions (even in the viscosity sense) are not clear. The approximation process is not easy as well and the numerical realization remains a challenge though this equation is a 
 scalar equation which is easier to discretize than the
vectorial one.  Here we used the discrete scheme described in \cite{Set} to get
numerical results.
\section{Results and discussion}
 
 \subsection{Explicit scheme}
 We briefly present the numerical scheme. We used an explicit scheme in time and the spatial  discretization has been performed  following  \cite{Set} .  
We set 
$$  \mathcal{G} = (FH)^* FH \mbox{ and } g^* = H^*F^*g~,$$
so that the equation (\ref{eq:levelset}) is 
$$
\frac{\partial\phi}{\partial t}=|\nabla \phi |\left( \lambda^2  \mathcal{G} \left(1_{\phi(t,\cdot)>0}\right)-\beta\mathrm{ div}\left(\frac{\nabla 
  \phi}{|\nabla \phi |} \right)- \lambda g^* \right).
$$
We set $t_n = n \Delta t, ~  \Phi^n  = \phi(t_n, \cdot),~X= (x_i,y_j)_{(i,j)\in I}$ with  $x_i = i \Delta x $ et $y_j = j \Delta y$.  The explicit Euler scheme gives :
$$\Phi^{n+1}(X) = \Phi^n(X) + \Delta t ~| \nabla \Phi^n|(  X)  \left ( \lambda ^2 \mathcal{G}(1_{\Phi^n >0}) -\beta \mbox{ curv}(\Phi^{n})  - \lambda g^*(t_{n}, X)\right )~.
$$
The curvature term is computed as 
$$ \mbox{curv }(\Phi) = \frac{ \Phi_{xx} (\Phi_y)^2 - 2 \Phi_{x} \Phi_y\Phi_{xy}  + \Phi_{yy} (\Phi_x)^2}{\left( (\Phi_x^2+ \Phi_y^2\right)^{3/2}} ~,$$
where $\Phi_x$  stands for the partial derivative with respect to $x$. 
  \noindent  The discrete approximation of the gradient is standard: 
 $$D^+_x\Phi (x,y) = \frac{\Phi(x+\Delta x , y)- \Phi(x,y)}{\Delta x},~D^-_x\Phi (x,y) = \frac{\Phi(x , y)- \Phi(x-\Delta x ,y)}{\Delta x}~;$$
$D^+_y \Phi$ are $D^-_y \Phi$ defined in the same way. A usual approximation for  \  $|\nabla \Phi| $  is given by : 
$$|\nabla \Phi|(X)\simeq \left[\max(D^+_x\Phi,0)^2+ \max(D^+_y\Phi,0)^2+\min(D^-_x\Phi,0)^2+\min(D^-_y\Phi,0)^2\right]^{1/2}~.$$
The non local term  $\mathcal{G}(1_{\Phi^n >0})$ is exactly computed.

\subsection{Numerical results}
The previous scheme has been implemented on a 3.6 GHz PC. A classical
reinitialization process has been used each 500 iterations. The test image size was  256 $\times$ 256 pixels. The other parameters of the computation were set to 
$$ \alpha = 10~,~\lambda = 2,~\Delta x = 1 \mbox{ and } \Delta t =
10^{-4}$$
and the blur kernel is a Gaussian kernel of standard deviation 5
pixels.

The computed image is quite satisfying (see figure 7.) However, we note a bad reconstruction along the symmetry axis due to   the  problem geometry  and a lack of information. Moreover we have to improve the algorithme behavior.  Indeed, we observe numerical instability (in spite of the regularization process) that leads to a very small time step choice. Therefore the computational time is quite long (about  2.5 hours). In addition, classical stopping criteria are not useful here~: the expected solution corresponds to a ``flat'' level of function $\Phi$ and the difference between two consecutive iterates means no sense. An estimate of the  cost function  decrease is not appropriate as well (we observe oscillations). We decided to stop after a  large enough number of iterations (here 20 000). 
 
In spite of all these disadvantages, this method is satisfactory considering the low signal to noise ratio of the radiograph. These good results can be explained by the strong assumptions
that we add (in particular the binary hypothesis) which are
verified by our synthetic object.  Anyway, the method has beeen successfully tested  on ``real''  images as well, that is imgaes of objets  with the same kind of properties (``almost'' binary)  but we cannot report them here (confidential data). 
A semi-implicit version of the algorithm is actually tested to improve stability.

\begin{figure}[H]
\begin{center}
\includegraphics[width=5cm,angle=270]{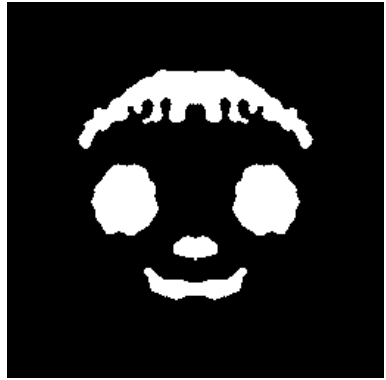}
\caption{Experimental results}\label{fig:results}
\end{center}
\end{figure}



\begin{thebibliography}{lll}

\bibitem{AA} {\bf I. Abraham, R. Abraham}
Technical Report CEA  (2001).
 
\bibitem{ambrosio}{\bf  L. Ambrosio,  N. Fusco et D. Pallara,}\textit{Functions of bounded variation and free discontinuity problems}, Oxford mathematical monographs, Oxford University Press, 2000.

\bibitem{BPDB} {\bf J.P. Bruandet, F. Peyrin, J.M. Dinten, M. Barlaud} {\sl
  3D tomographic reconstruction of binary images from cone-beam
  projections: a fast level-set approach}\\
2002 IEEE International Symposium on Biomedical Imaging, p. 677-80, (2002)

\bibitem{CKP}{\bf  E. Casas, K. Kunisch and C. Pola},\textit{ Regularization by Functions of Bounded Variation and
Applications to Image Enhancement}, Applied  Mathematics and  Optimization,  40:229Ð257 (1999)



\bibitem{Din} {\bf J.-M.  Dinten} {\sl Tomographie \`a partir d'un nombre limit\'e de projections
: r\'egularisation par champs markovien}\\
PHD thesis, Universit\'e d'Orsay Paris-Sud (1990)

\bibitem{Dus} {\bf N. J. Dusaussoy} {\sl Image reconstruction from
  projections}\\
SPIE's international symposium on optics, imaging and instrumentation. San Diego
(1994)

\bibitem{FKC}{\bf  H. Feng,  W. Karl,  D. Castanon} {\sl A curve
  evolution approach to object-based tomographic reconstruction}\\
IEEE Trans. on Image Proc., 12, 44-57, (2003)

\bibitem{Han}{\bf K. Hanson} {\sl Tomographic reconstruction of
  axially symmetric objects from a single radiograph}\\
High Speed Photography, 491, (1984)

\bibitem{Her}{\bf G. Herman} {\sl Image reconstruction from projections: the fundamentals of
computerized tomography}\\
Academic Press (1980)

\bibitem{Lag} {\bf J.M. Lagrange} {\sl Reconstruction tomographique \`a partir d'un petit nombre de vues}\\
PHD thesis, ENS Cachan (1998)

\bibitem{MS} {\bf D. Mumford - J. Shah} {\sl Optimal approximations by
  piecewise smooth functions and associated variational problem}\\
Comm. Pure and Appl. Math. 42, 577-685 (1989)

\bibitem{Set} {\bf J.A. Sethian} {\sl Theory, algorithm and
  applications of level set method for propagating interfaces}\\
Iserles, A. (ed.), Acta Numerica Vol. 5, 1996. Cambridge: Cambridge
  University Press. 309-395 (1996)
\bibitem{Stein} {\bf E. Stein} {\sl  
 Singular integrals and differentiability properties of functions},  Princeton University Press  
\end{thebibliography}
\end{document}